\mag=\magstep 1     

\documentclass{article}  
\usepackage{amsmath,amssymb, mathrsfs}
\usepackage{romp}

\title{On the sensitivities dependence in non-autonomous dynamical systems }

\begin{document}

\maketitle

\hspace*{0.6cm} \begin{center}{\large Yang Chengyu}$^{\: a}$\symbolfootnote[1]{Corresponding author.
tel.: (+86) 1389*******.\\
\hspace*{0.25cm} {\it E-mail addresses:}
y582365724@126.com (C. Yang), lizm@nwu.edu.cn (Z. Li).\\
\hspace*{0.4cm}$^1$Supported in Part by the National Natural Science Foundation of China (NSFC No. 11301417).}, {\large Li Zhiming}$^{\: b,1}$\end{center}

\bigskip

\hspace*{0.4cm} $^a$ {\it \small School of Economic and Management, Northwest University, Xian, Shaanxi 710127, China}\\
\hspace*{0.4cm} $^b$ {\it \small Department of Mathematics, Northwest University, Xian, Shaanxi 710127, China}\\

\bigskip
\begin{abstract}
For discrete autonomous dynamical systems (ADS) $(X, d, f)$, it was found that in the three conditions defining Devaney chaos, topological transitivity and dense periodic points together imply sensitive dependence on initial condition(Banks, Brooks, Cairns, Davis and Stacey, 1992). In this paper, the result of Banks et al. is generalized to a class of the non-autonomous dynamical systems (NADS) $(X,f_{1,\infty})$. Also, by the studying of NADS over their iterated systems $(X,f_{1,\infty}^{[k]})$, we know that for two sensitive NADS, the one which preserve sensitive in its any times iterated systems is more sensitive than the one not. In this case, several sufficient conditions ensuring two kinds of sensitivities are preserved under the arbitrary number of iterations of certain NADS are given.

\end{abstract}

\bigskip
\noindent {\it keywords:} Non-autonomous dynamical system; topological transitivity; sensitive dependence on initial condition; collective sensitivitive; 
Banks Theorem

\section{\hspace*{-0.15in}. Introduction}
In 1971, Ruelle introduced the first precise definition for sensitivity \cite{Ruelle}. In 1986, Devaney proposed the widely accepted definition of chaos (topological transitivity, dense periodic points and sensitivity), and emphasized the significance of sensitivity in describing ADS \cite{Devaney}. Afterwards, Li-Yorke sensitivity \cite{Akin}, $n$ sensitivity \cite{Xiong}, and collective sensitivity \cite{Wang} were
successively proposed, and each of these concepts is used to describe the complexity of dynamical systems. Moreover, in 1992, Banks et al found that in the three conditions defining Devaney chaos in the ADS,
topological transitivity and dense periodic points together imply sensitivity \cite{Banks}.

To NADS, in 1996, Kolyada and Snoha\cite{Kolyada} investigated the topological entropy and properties
for NADS. In 2011, C\'{a}novas introduced the definition of chaos in NADS \cite{C1}, and
studied the relationship between the chaos and the topological entropy.
Dvor\v{a}kova studied the relation between NADS $(X,f_{1,\infty})$
and ADS $(X, f)$, where all the involved maps are defined
on the closed unit interval $[0, 1]$  and $f_{1, \infty}$ converges uniformly to
$f$ \cite{Dvorakova}. In 2012, Balibrea and Oprocha explored the properties of Li-Yorke chaos in NADS and studied the
relation between topological weak mixing and topological entropy \cite{Balibrea}. In 2013, Murillo-Arcila studied
the topological mixing for linear NADS, and proved that the $n$th topological mixing does not
imply the $(n+1)$th, which is different from the corresponding result in ADS \cite{Murillo}.

Regarding the sensitivity in NADS, some results were obtained. In 2006, Tian and Chen introduced the definition
of sensitivity for NADS \cite{Tian}. Then, \cite{Shi}
proposed the concept of Devaney chaos for NADS, and asked as an open problem that whether or not the
previously stated theorem by Banks et al. can be generalized from DAS to NADS. In 2013, Wu and Zhu \cite{Wu} studied the
relative hereditary property of sensitivity in NADS defined on compact metric spaces, and proved that, when the mapping sequence $(X, f_{1, \infty})$
converges uniformly, for any positive integer $k$, the system $(X, f_{1, \infty}^{[k]})$ (definition given in Section 2) also holds the sensitivity.
In addition, Wu and Zhu found a sufficient condition ensuring
relative hereditary property of sensitivity.

In this paper, some necessary definitions like the concepts of collective sensitivity dependence on initial condition (abbreviated as collective sensitivity) and synchronous sensitivity dependence on initial condition (abbreviated as
synchronous sensitivity) are introduced in section 2. And in section 3, two kinds of issues are researched:

One hand, for a NADS $(X, f_{1,\infty})$ defined on a metric space (no matter compact or not) and having sensitivity dependence on initial condition(or collective sensitivity), if its mapping sequence $f_{1, \infty}$ is either finitely
generated(in definition \ref{f1}) or converges uniformly to a map $f$, to every positive integer $k$, the system $(X, f_{1, \infty}^{[k]})$ hold the sensitivity.

In the other hand, the result of Banks et al is generalized to the finitely generated NADS. In addition, for linear NADS, we also investigate the relationship between topological transitivity and sensitivities.

\section{\hspace*{-0.15in}. Preliminaries}

 In this section, we mainly recall some relevant concepts and terminologies.

\begin{definition}\cite{Kolyada}
Let $(X,d)$ be a metric space and
$f_{i}: X\rightarrow X$, $i = 1, 2, \cdots$, a sequence of continuous maps.
We write $\{f_{i}\}_{1}^{\infty}=f_{1, \infty}$.
For any $x\in X$, the the following sequence
                                $$x, f_{1}(x), f_{1}^{2}(x), f_{1}^{3}(x), \cdots, f_{1}^{n}(x), \cdots,$$
is denoted by $tra(x)$. Here $f_{i}^{n}=f_{n-1+i}\circ f_{n-2+i}\circ \cdots \circ f_{i+1}\circ f_{i}$,
$f_{i}\in f_{1, \infty}$ and $f_{i}^{0}=id$, $i=1, 2, ..., n$. We call
$ f_{1, \infty}$ a {\bf NADS} on $X$ and denote this system by $(X, f_{1, \infty})$.
The set of points determined by {\bf $tra(x)$} is the orbit of $x$ which denoted by {\bf $orb(x)$}.

\end{definition}

\begin{definition}\cite{Murillo}
Let $(X,f_{1,\infty})$ be a NADS. A subset $A$ of $X$ is said to be an {\bf invariant set} for
$(X, f_{1,\infty})$ if $f_{i}(A)\subset A$ for every $f_i \in f_{1, \infty}$.
\end{definition}

\begin{definition}
 A system $(X, f_{1, \infty})$ is periodic if there is a positive integer $p$ such that, for any positive integer $k$, it holds that $f_{k+p}=f_{k}$. The least such number $p$ is called the period of $(X, f_{1, \infty})$.
Specially, if $(X,f_{1, \infty})$ is periodic with period 1, it is an ADS.
\end{definition}

\begin{definition}\label{f1}
$(X,f_{1, \infty})$ said to be a {\bf finitely generated} NADS if there exists a finite set $F$ of continuous maps on $X$ such that, every $f_i$ of $f_{1, \infty}$ belongs to $F$.
\end{definition}

All periodic NADS are finitely generated, but not vice versa.

\begin{definition}
$(X,f_{1,\infty})$ is said to be {\bf commutative} if for any pair of positive integer $m,n$ and any $x\in X$, $f_{1}^{n}\circ f_{1}^{m}(x)=f_{1}^{m}\circ f_{1}^{n}(x)$.
\end{definition}

\begin{definition}\cite{C2}
A point $x$ of $X$ is a {\bf periodic point} of $f_{1, \infty}$ if there is a positive integer N such that, for any natural number $k$, it holds that $f_{1}^{k}(x)=f_{1}^{N+k}(x)$. The least such number N is called the period of $x$. If the period is 1, it is a fixed point. The set of periodic points is denoted by $P(f_{1, \infty})$.
\end{definition}

In any ADS, the orbit of any periodic point forms an invariant set. However, this no longer holds for all NADS. A counterexample is provided below.

\begin{example}\label{example}
Let $X=\mathbb{R}$, and four continuous maps $f_1, f_2, g_1$ and $g_2$ on $X$ are defined as follows: $f_{1}(x)=x+1, f_{2}(x)=x-1, g_{1}(x)=x^{2}$, and $g_{2}(x)=-x$.
Then we get two sequences of maps: $f_{1, \infty}=(f_{1}, f_{2}, f_{1}, f_{2}, \cdots)$
and $g_{1, \infty}=(g_{1}, g_{2}, g_{1}, g_{2}, \cdots)$. In this case, $(X,f_{1, \infty})$ and $(X,g_{1, \infty})$ are two periodic NADS on $X$.

In $(X,f_{1, \infty})$, for any $x \in X$, we have $tra(x)=(x, x+1, x, x+1, x, x+1, \cdots)$. Hence it holds that $x\in P(f_{1, \infty})$.
Particularly, $tra(1) =(1, 2, 1, 2, \cdots), orb(1)=\{1, 2\}$, $f_{1}(orb(1))=\{2, 3\}\not\subseteq orb(1)$, $f_{2}(orb(1))=\{0,1\}\not\subseteq orb(1)$.
So, $orb(1)=\{1,2\}$ is not an invariant set of $f_{1, \infty}$.

In $(X, g_{1, \infty})$, $tra(-1)=(-1, 1, -1, 1, \cdots), orb(1)=\{-1, 1\}$. This means that -1 is a periodic point of $g_{1, \infty}$ with period of 2.
Furthermore, it also satisfies
$ g_{1}(orb(-1))=\{1, -1\} \subseteq orb(-1), g_{2}(orb(-1))=\{1, -1\}\subseteq orb(-1)$. Hence, $orb(-1)$ is  an invariant set for $g_{1, \infty}$.
\end{example}

 The above examples demonstrate a fact that properties of the set of periodic points in NADS are quite different from what in an ADS.  Consequently, it is necessary to define a new kind of periodic point with a stronger property.

\begin{definition}
We call $x$ an {\bf invariant periodic point} of $f_{1, \infty}$
if $x\in P(f_{1, \infty})$ and $orb(x)$ is an invariant set for $f_{1, \infty}$.
\end{definition}

\begin{definition}
A system $(X, f_{1,\infty})$ is {\bf transitive} if for any pair $U$ and $V$ of nonempty open subsets of $X$, there exists $N>0$
so that $f_{1}^{N}(U)\bigcap V\not=\emptyset$.
\end{definition}

\begin{definition}\cite{Grosse}
Let $X$ be a separable Frechet space, $d$ a translation invariant metric on X, and $f_{1, \infty}$  a sequence of continuous maps on X.
Then we call $(X, d, f_{1, \infty})$ a {\bf linear NADS}.
\end{definition}

\begin{definition}\cite{Shi}
$(X,f_{1,\infty})$ is said to have {\bf sensitivity dependence on initial conditions} if there exists some $\delta>0$ such that, for every $x\in X$ and $\epsilon>0$, there exists $y\in X$ with $d(x,y)<\epsilon$ such that, for some positive integer $n$, $d(f_{1}^{n}(x),f_{1}^{n}(y))>\delta$. The number $\delta$ is called a sensitivity constant for $(X,f_{1,\infty})$.

\end{definition}

\begin{definition}\cite{Wang}
$(X,f_{1,\infty})$ is said to have {\bf collective sensitivity} if there exists some $\delta>0$ such that,for any finitely many distinct points $x_{1},x_{2},x_{3},\cdots,x_{n}$ of $X$ and a arbitrary $\epsilon>0$,there exists same number of distinct points $y_{1},y_{2},y_{3},\cdots,y_{n}$ of $X$ and some positive integer $k$ satisfying the following two conditions£º

(1)$d(x_{i},y_{i})<\epsilon,i=1,2,3,\cdots,n;$

(2)there exists an $i_{0}$ with $1\leq i_{0}\leq n$ such that $ d(f_{1}^{k}(x_{i}),f_{1}^{k}(y_{i_{0}}))\geq \delta,i=1,2,3,\cdots,n,$ or $ d(f_{1}^{k}(y_{i}),f_{1}^{k}(x_{i_{0}}))\geq \delta,i=1,2,3,\cdots,n$.\\
The number $\delta$ is called a collective sensitivity constant for $(X,f_{1,\infty})$.
\end{definition}

The Theorem 2.3.in paper\cite{Wang} proved that an ADS having collective sensitive and its induced (sub)hyperspace dynamical systems equipped with the Vietoris topology having sensitive dependence on initial conditions are equivalent conditions. This results also holds in NADS\cite{yang}.

\begin{definition}
$(X,f_{1,\infty})$ is said to have {\bf synchronous sensitivity} if there exists some $\delta>0$ such that,for any finitely many distinct points $x_{1},x_{2},x_{3},\cdots,x_{n}$ of $X$ and a arbitrary $\epsilon>0$,there exists same number of distinct points $y_{1},y_{2},y_{3},\cdots,y_{n}$ of $X$ and some positive integer $k$ satisfying the following two conditions

(1)$d(x_{i},y_{i})<\epsilon,i=1,2,3,\cdots,n;$

(2)$ d(f_{1}^{k}(x_{i}),f_{1}^{k}(y_{i}))\geq \delta,i=1,2,3,\cdots,n$.\\
The number $\delta$ is called a synchronous sensitivity constant for $(X,f_{1,\infty})$
\end{definition}

\begin{definition}\cite{Shi}
$(X,f_{1,\infty})$ is said to be {\bf Devaney chaos}, if it satisfies the following three conditions:

(1)$\overline{P(f_{1,\infty})}=X$

(2)$f_{1,\infty}$ is topologically transitive

(3)$f_{1,\infty}$ has sensitive dependence on initial conditions.

\end{definition}

\begin{definition}\cite{Wu}
Let $(X,f_{1,\infty})$ be a NDAS, for any positive integer $k$, denote $$(f_{1}^{k},f_{k+1}^{k},\cdots,f_{km+1}^{k},\cdots)
=f_{1,\infty}^{[k]},$$ we know
$(X,f_{1,\infty}^{[k]})$ also is a NDAS, and we call this {\bf $k$-th iterate system} of $(X,f_{1,\infty})$.
\end{definition}

For an ordinary NDAS, it can have infinite numbers of iterate systems.

\section{\hspace*{-0.15in}. Main Results}

Just from the definition of sensitivity, we can find out that in a NADS $(X,f_{1,\infty})$, if one of its iterate system is sensitive, it is sensitive too, but which is not vice versa. In paper \cite{Wu}, Wu show us that there is a sensitive NADS $(Y,g_{1,\infty})$, for any positive integer $m$, its iterate system $(Y,g_{1,\infty}^{[m]})$ is not sensitive. So now, we can distinguish two sensitive NADS which is more chaotic by studying whether its iterate system is sensitive.

What is more, in following theorems from \ref{th3.1} to \ref{th3.99}, several sufficient conditions ensuring that sensitivity is preserved under any iterate system of a certain NADS will be given.

\bigskip
\begin{theorem}\label{th3.1}
Let $(X,f_{1,\infty})$ be a a finitely generated NADS and for any $ f_{i}\in f_{1,\infty}$,$f_{i}$ is uniform continuity. If $f_{1,\infty}$ has sensitivity dependence on initial conditions, then for every positive integer $k$, so does $f_{1,\infty}^{[k]}$.
\end{theorem}

{\bf Proof:} Let $\delta>0$ be a sensitivity constant for $(X,f_{1,\infty})$ and take any given $k>1$.

As for any $ f_{i}\in f_{1,\infty}$, $f_{i}$ is uniform continuity, we obtain that for every positive integer i, $0\leq l\leq k+2$, $f_{i}^{l}$ is uniform continuity.

Because $f_{1,\infty}$ is finitely generated, there exists $\epsilon_{\delta}>0$ which ensure that for any $x, y\in X, i>0$, once $d(x,y)<\epsilon_{\delta},0\leq l\leq k+2$, we always have $$d(f_{1}^{l}(x),f_{1}^{l}(y))<\delta.\eqno(3.1)$$

Now, we will prove that $\epsilon_{\delta}$ is a sensitivity constant for $(X,f_{1,\infty}^{[k]})$. For any $x\in X$ and $\epsilon>0$ (might take $\epsilon<\frac{\epsilon_{\delta}}{4}$ here),
as $(X,f_{1,\infty})$ is sensitive, there are positive integer $n_{x,\epsilon}$ and $y\in X$ satisfying $$d(x,y)<\epsilon,d(f_{1}^{n_{x,\epsilon}}(x),f_{1}^{n_{x,\epsilon}}(y))>\delta$$ ( $n_{x,\epsilon}$ is the smallest positive integer number to satisfy this inequality here).
By (3.1), we know $n_{x,\epsilon}>k+2$.

Because of $n_{x,\epsilon}>k+2>k$,
we can find a positive integer $r$ satisfying $n_{x,\epsilon}=rk+q,0\leq q\leq k-1$.

As
$$\delta<d(f_{1}^{n_{x,\epsilon}}(x),f_{1}^{n_{x,\epsilon}}(y))=d(f_{rk+1}^{q}(f_{1}^{rk}(x)),f_{rk+1}^{q}(f_{1}^{rk}(y)))$$
and $q=n_{x,\epsilon}-rk\leq k-1<k+2$,
by considering the inequality (3.1), $d(f_{1}^{rk}(x),f_{1}^{rk}(y))>\epsilon_{\delta}$ is obvious.

Because of the arbitrariness of $x$ and $\epsilon$,
$\epsilon_{\delta}$ is a sensitivity constant to $(X,f_{1,\infty}^{[k]})$, and for the arbitrariness of $k$, theorem holds.$\Box$

\bigskip

\begin{corollary}
Let $(X,f_{1,\infty})$ be a finitely generated NADS on compact space $X$. If $f_{1,\infty}$ have sensitivity dependence on initial conditions, then for every positive integer $k$, so does $f_{1,\infty}^{[k]}$.
\end{corollary}

\bigskip

To collective sensitive NDAS, we have some similar conclusions here too.

\begin{theorem}\label{th3.2}
Let $(X,f_{1,\infty})$ be a finitely generated NADS and for any $ f_{i}\in f_{1,\infty}$, $f_{i}$ is uniform continuity. If $f_{1,\infty}$ has collective sensitivity, then for every positive integer $k$, so does $f_{1,\infty}^{[k]}$.
\end{theorem}

{\bf Proof:}
Let $\delta>0$ be a collective sensitivity constant of $f_{1,\infty}$, and take any given $k>1$.

Same as theorem\ref{th3.1}, for the given $\delta$ and $k$, there is $\epsilon_{\delta}>0$ satisfying that
for any $x,y\in X,i>0$, once $d(x,y)<\epsilon_{\delta},0\leq l\leq k+2$, we always have $$d(f_{1}^{l}(x),f_{1}^{l}(y))<\delta.\eqno(3.2)$$

Now, we will prove that $\epsilon_{\delta}$ is a collective sensitivity constant of $f_{1,\infty}$.

As $f_{1,\infty}$ have collective sensitivity, for any finitely many distinct points $x_{1},x_{2},x_{3},\cdots,x_{n}$ of $X$ and a arbitrary $\epsilon>0$ ($\epsilon<\epsilon_{\delta}$), there exists same number of distinct points $y_{1},y_{2},y_{3},\cdots,y_{n}$ of $X$ and some positive integer $n_{\epsilon}$ ensuring:

(1)$d(y_{i},x_{i})<\epsilon,i=1,2,3,\cdots,n $

(2)there exists $1\leq i_{0}\leq n$ such that, $ d(f_{1}^{n_{\epsilon}}(x_{i}),f_{1}^{n_{\epsilon}}(y_{i_{0}}))\geq  \delta,i$=$1,2,3,\cdots,n$; or $ d(f_{1}^{n_{\epsilon}}(y_{i}),f_{1}^{n_{\epsilon}}(x_{i_{0}}))\geq \delta,i$=$1,2,3,\cdots,n$.

As $\epsilon<\epsilon_{\delta}$, by considering $(3.2)$, we have $n_{\epsilon}>k+2>k$. So there is a positive integer $r$ ensuring $n_{\epsilon}=rk+q, 0\leq q\leq k-1.$

When $d(f_{1}^{n_{\epsilon}}(x_{i}),f_{1}^{n_{\epsilon}}(y_{i_{0}}))\geq \delta,i=1,2,3,\cdots,n$ established, we have
\begin{eqnarray*}
\ d(f_{1}^{n_{\epsilon}}(x_{i}),f_{1}^{n_{\epsilon}}(y_{i_{0}}))=d(f_{1}^{rk+q}(x_{i}),f_{1}^{rk+q}(y_{i_{0}}))
=d(f_{rk+1}^{q}(f_{1}^{rk}(x_{i})),f_{rk+1}^{q}(f_{1}^{rk}(y_{i_{0}})))\geq \delta,i=1,2,3,\cdots,n
\end{eqnarray*}

Since $q$=$n_{x,\epsilon}-rk\leq k-1<k+2$,and $(3.2)$, we have $d(f_{1}^{rk}(x_{i}),f_{1}^{rk}(y_{i_{0}}))>\epsilon_{\delta},i$=$1,2,3,\cdots,n.$

In the same way, when $ d(f_{1}^{n_{\epsilon}}(y_{i}),f_{1}^{n_{\epsilon}}(x_{i_{0}}))\geq \delta,i$=$1,2,3,\cdots,n$ established, we have
$$d(f_{1}^{rk}(x_{i_{0}}),f_{1}^{rk}(y_{i}))>\epsilon_{\delta},i=1,2,3,\cdots,n.$$

Because of the arbitrariness of $x_{1},x_{2},x_{3},\cdots,x_{n}$
and $\epsilon$, we know $\epsilon_{\delta}$ is a collective sensitivity constant of $(X,f_{1,\infty}^{[k]})$. And for the arbitrariness of $k$, theorem holds.$\Box$

\bigskip

\begin{corollary}
Let $(X,f_{1,\infty})$ be a finitely generated NADS on compact space $X$. If $f_{1,\infty}$ have collective sensitivity, then for every positive integer $k$, so does $f_{1,\infty}^{[k]}$.
\end{corollary}

Let $f_{1,\infty}$ be a sequence of continuous maps on metric space $X$, $f$ is a map on $X$, we said $f_{1,\infty}$ converges uniformly to $f$, if for any $\epsilon>0$, there exists a positive integer $N_{0}$, satisfying for any $x\in X$ and $n>N_{0}$, we all have $d(f_{n}(x),f(x))<\epsilon$.

\begin{theorem}\label{th3.3}
Assume $f_{1,\infty}$ converges uniformly to $f$ and for any $ f_{i}\in f_{1,\infty}$, $f_{i}$ is continues uniformly. Then:

(1) For any positive integer $k$, $\{f_{n}^{k}\}_{n=1}^{\infty}$ converges uniformly to $f^{k}$;

(2) For any $\epsilon>0$ and positive integer $k$, there exists $\delta(\epsilon)>0$ and positive integer $N(k)$ such that for any pair $x,y\in X$ with $d(x,y)<\delta(\epsilon)$ and any $n\geq N(k)$, $d(f_{n}^{k}(x),f_{n}^{k}(y))<\frac{\epsilon}{2}$.
\end{theorem}

{\bf Proof:}(1)$k=1$, $f_{1,\infty}$ converges uniformly to $f$; $k=2$, for any $\epsilon>0$, any $x\in X$,
\begin{eqnarray*}
\ d(f_{n}^{2}(x),f^{2}(x))&\leq &d(f_{n}^{2}(x),f(f_{n}(x)))+d(f(f_{n}(x)),f^{2}(x))\\
&=&d(f_{n+1}(f_{n}(x)),f(f_{n}(x))+d(f(f_{n}(x)),f^{2}(x)).
\end{eqnarray*}

As $f_{1,\infty}$ converges uniformly to $f$,
there exists a positive integer $N_{0}$ satisfying for any $X\in X$ and $n>N_{0}$, we have $d(f_{n}(x),f(x))<\frac{\epsilon}{2}$. Especially, $d(f_{n+1}(f_{n}(x)),f(f_{n}(x)))<\frac{\epsilon}{2}$.

As $ f_{i}\in f_{1,\infty}$, $f_{i}$ is continues uniformly, there exists $\delta>0$ satisfying that as long as $d(f_{n}(x),f(x))<\delta$, we  have $ d(f(f_{n}(X)),f^{2}(x))<\frac{\epsilon}{2}$.

Since $f_{1,\infty}$ converges uniformly to $f$, there also exists a positive integer $N_{1}$ satisfying that for any $x\in X$ and $n>N_{1}$, we all have
$d(f_{n}(x),f(x))<\delta$. Then $d(f(f_{n}(x)),f^{2}(x))<\frac{\epsilon}{2}.$

In summary, there is $N=max\{N_{0},N_{1}\}$ such that, for any $x\in X$ and positive integer $n>N$, we have
\begin{eqnarray*}
\ d(f_{n}^{2}(x),f^{2}(x))&\leq &d(f_{n}^{2}(x),f(f_{n}(x)))+d(f(f_{n}(x)),f^{2}(x))\\
&=&d(f_{n+1}(f_{n}(x)),f(f_{n}(x))+d(f(f_{n}(x)),f^{2}(x))\leq \frac{\epsilon}{2}+\frac{\epsilon}{2}=\varepsilon.
\end{eqnarray*}

It means that $\{f_{n}^{2}\}_{n=1}^{\infty}$ converges uniformly to $f^{2}$.

Assume for positive integer $k>1$, $\{f_{n}^{k}\}_{n=1}^{\infty}$ converges uniformly to $f^{k}$; we can prove that $\{f_{n}^{k+1}\}_{n=1}^{\infty}$ converges uniformly to $f^{k+1}$ in the same way.

(2)for any positive integer $k$, $\epsilon>0,x,y\in X$,
\begin{eqnarray*}
\ d(f_{n}^{k}(x),f^{k}_{n}(y))\leq d(f_{n}^{k}(x),f^{k}(x)))+d(f^{k}(x),f^{k}(y))+d(f^{k}(y)),f^{k}_{n}(y)).
\end{eqnarray*}

Since $\{f_{n}^{k}\}_{n=1}^{\infty}$ converges uniformly to $f^{k}$, there exists positive integer $N(k)$ satisfying for any $x\in X$ and $n>N(k)$, we have $d(f_{n}^{k}(x),f^{k}(x))<\frac{\epsilon}{6}$.

For the uniform continuity, there also exists $\delta(\epsilon)>0$, as long as $d(x,y)<\delta(\epsilon)$, we have $d(f^{k}(x),f^{k}(y))<\frac{\epsilon}{6}$.

Then there exists $\delta(\epsilon)>0$ and positive integer $N(k)$ such that, for any $x,y\in X$ with $d(x,y)<\delta(\epsilon)$ and any $n>N(k)$, we have
\begin{eqnarray*}
\ d(f_{n}^{k}(x),f^{k}_{n}(y))&\leq &d(f_{n}^{k}(x),f^{k}(x)))+d(f^{k}(x),f^{k}(y))+d(f^{k}(y)),f^{k}_{n}(y))\\
&\leq & \frac{\epsilon}{6}+\frac{\epsilon}{6}+\frac{\epsilon}{6}=\frac{\epsilon}{2}.
\end{eqnarray*}$\Box$

Especially, if $X$ is compact, theorem \ref{th3.3} is lemma 2.1 in paper \cite{Wu}.

\bigskip

\begin{theorem}\label{th3.4}
Let $(X,f_{1,\infty})$ be a NADS which sequence of maps converges uniformly to continues map $f$. And for any $ f_{i}\in f_{1,\infty}$, $f_{i}$ is continues uniformly.
If $f_{1,\infty}$ has collective sensitivity, then for every positive integer $k$, so does $f_{1,\infty}^{[k]}$.
\end{theorem}

{\bf Proof:}
Assume $\delta>0$ is a collective sensitivity constant for $f_{1,\infty}$. Taking any integer number $k>1$, according to theorem\ref{th3.3}(2), for $2\delta$ and $k$, there exists $\epsilon_{\delta}>0$ and positive integer $n_{0}$($n_{0}>3k$) ensure that, for any $x,y\in X,d(x,y)<\varepsilon_{\delta}$ and any $n\geq n_{0}$, we have
$$d(f_{n}^{i}(x),f^{i}_{n}(y))<\delta,i=1,2,3,\cdots,k\eqno (3.3)$$

Now, we will prove that $\epsilon_{\delta}$ is a collective sensitivity constant of $(X,f_{1,\infty}^{[k]})$.

For any $0<i\leq 2n_{0}, f_1^{i}$ is continues uniformly. Then there exists $\epsilon^{*}>0$ satisfying for any $0<i\leq 2n_{0}$ and any $x,y\in X$ with $d(x,y)<\epsilon^{*}$, we have $d(f_{1}^{i}(x),f_{1}^{i}(y))<\delta, 0<i\leq 2n_{0}$.
Since $f_{1,\infty}$ have collective sensitive, for any finitely many distinct points $x_{1},x_{2},x_{3},\cdots,x_{n}$ of $X$ and a arbitrary $\epsilon>0$ ($\epsilon<\epsilon^{*}$),
there exists the same number of distinct points $y_{1},y_{2},y_{3},\cdots ,y_{n} \in X$ and $m>0$ such that

(1)$d(y_{i},x_{i})<\epsilon,i=1,2,3,\cdots,n;$

(2)there exists $1\leq i_{0}\leq n$ making such that $ d(f_{1}^{m}(x_{i}),f_{1}^{m}(y_{i_{0}}))\geq \delta,i=1,2,3,\cdots,n;$ or $ d(f_{1}^{m}(y_{i}),f_{1}^{m}(x_{i_{0}}))\geq \delta,i=1,2,3,\cdots,n$.\\
By considering the choice of $\epsilon$, we know $m>2n_{0}>6k$. Then there exists positive integer $p\geq 6$ satisfying
$m=pk+q,1\leq q\leq k-1$.

Since $m>2n_{0}>6k$ and $1\leq q\leq k-1$, we get $pk+1>n_{0}$.

When $ d(f_{1}^{m}(x_{i}),f_{1}^{m}(y_{i_{0}}))\geq \delta,i=1,2,3,\cdots,n$ established, we botain
\begin{eqnarray*}
\ d(f_{1}^{m}(x_{i}),f_{1}^{m}(y_{i_{0}}))&=&d(f_{1}^{pk+q}(x_{i}),f_{1}^{pk+q}(y_{i_{0}}))\\
&=&d(f_{pk+1}^{q}(f_{1}^{pk}(x_{i})),f_{pk+1}^{q}(f_{1}^{pk}(y_{i_{0}})))\geq \delta,i=1,2,3,\cdots,n.
\end{eqnarray*}

By (3.3), as $pk+1>n_{0},q<k$, $d(f_{1}^{kp}(x_{i},f_{1}^{kp}(y_{i_{0}})\geq \epsilon_{\delta},i=1,2,3,\cdots,n.$\\

In the same way, when $ d(f_{1}^{m}(y_{i}),f_{1}^{m}(x_{i_{0}}))\geq \delta,i=1,2,3,\cdots,n$ established, we have
$$ d(f_{1}^{m}(y_{i}),f_{1}^{m}(x_{i_{0}}))\geq \epsilon_{\delta},i=1,2,3,\cdots,n.$$

Because of the arbitrariness of $x_{1},x_{2},x_{3},\cdots,x_{n}$ and $\epsilon$,
$\epsilon_{\delta}$ is a collective sensitivity constant of $(X,f_{1,\infty}^{[k]})$.$\Box$

\bigskip

\begin{corollary}\label{th3.99}
Let $(X,f_{1,\infty})$ be a NADS which sequence of maps converges uniformly to continues map $f$ and $X$ is compact.
If $f_{1,\infty}$ has collective sensitivity, then for every positive integer $k$, so does $f_{1,\infty}^{[k]}$.
\end{corollary}

\bigskip

Especially, when a NADS degenerates into a ADS, all conclusions above are available.

\bigskip

\begin{theorem}\label{th3.6}\cite{Banks}
$(X,f)$ is a autonomous dynamical system without isolate point which satisfies:\\
(1)Periodic points are dense in $X$($\overline{P(f)}=X$);\\
(2)$(X,f)$ is topological transitive.\\
Then, $(X,f)$ have sensitivity dependence on initial conditions.
\end{theorem}

Theorem \ref{th3.6} is the Banks-Brooks-Cairns-Davis-Stacey theorem which is not only simplifying the definition of Devaney chaos but also showing us the relationship between transitivity, periodic points and sensitivity in ADS. In \cite{Shi}, we know that the Banks et al theorem holding or not in NADS is still an open question. In the following theorem of this paper, we will prove it in a class of NADS firstly.

\begin{theorem}\label{th3.7}
Let $X$ be a metric space without isolate point.If a finitely generated non-autonomous dynamical system $(X,f_{1,\infty})$ satisfies the following conditions :\\
 (1)$f_{1,\infty}$is topologically transitive,\\
 (2)$\overline{P(f_{1,\infty})}=X$ ,\\
 (3)existing two invariant periodic points $x,y \in X$ and $orb(x)\cap orb(y)=\emptyset$
 \\
then the system is sensitive .
\end{theorem}

{\bf Proof:} For any nonempty subset $A$ of $X$ and $x\in X$, write $d(x,A)=inf\{d(x,y)|y\in X\}$ (If $A$ is finite, then $d(x,A)=min\{d(x,y)|y\in X\}$).

Let the two invariant periodic points is $p_{1}, p_{2}$ and $orb(p_{1})\cap orb(p_{2})=\emptyset$. Then we can note $$\delta=\frac{1}{3}min\{d(x,y)|x\in orb(p_{1}),y\in orb(p_{2})\}>0.$$

Firstly, we will prove that for any $x\in X$, there is $p_{x}\in \{p_{1},p_{2}\}$ satisfying $d(x,orb(p_{x}))>\delta$.

Three cases could happen here:

(1)$x\in orb(p_{1})$,$d(x,orb(p_{2})\geq 3\delta>\delta$.

(2)$x\in orb(p_{2})$,$d(x,orb(p_{1})\geq 3\delta>\delta$.

(3)$x\not\in orb(p_{1})\bigcup orb(p_{2})$, there is $a\in orb(p_{1}), b\in orb(p_{2})$ satisfying
$$d(x,a)=d(x,orb(p_{1})), d(x,b)=d(x,orb(p_{2})), d(a,b)\geq 3\delta.$$
For the triangle inequality:$$d(x,a)+d(x,b)\geq d(a,b)\geq 3\delta,$$
we have $d(x,a)\geq\frac{3\delta}{2}>\delta$ or $d(x.b)\geq\frac{3\delta}{2}>\delta$.

Because of the selection of $a,b$, we know $d(x,orb(p_{1}))>\delta$ or $d(x,orb(p_{2}))>\delta$.

In this case, for any $x\in X$, there is $p_{x}\in \{p_{1},p_{2}\}$ satisfying $d(x,orb(p_{x}))>\delta$.

Furthermore, as $p_{1}$,$p_{2}$ are invariant periodic points,
which means for any $f_{i}\in f_{1,\infty}$, we have
$$f_{i}(orb(p_{t}))\subset orb(p_{t}),t=1,2.$$

Then for any positive integer $i,m$,
$$d(x,f_{i}^{m}(orb(p_{x})))\geq d(x,orb(p_{x}))>\delta.$$

Noting $\eta=\frac{\delta}{8}$, we now will prove that $\eta$ is a sensitivity constant to $(X,f_{1,\infty})$.

For any $\epsilon>0$(let $\epsilon<\eta$ here), $\overline{P(f_{1,\infty})}=X$ and X has no isolated point, so there is  $q\in P(f_{1,\infty})$ satisfying
$q\in S(x,\epsilon),q\not=x$.

Here $S(x,\epsilon)=\{y\in X|d(x,y)<\epsilon\}$, and let $q$ has period $n_{q}$,

$$V(p_{x})=\bigcap_{i>0,n_{q}\geq l\geq 0} (f_{i}^{l})^{-1}(S(f_{i}^{l}(p_{x}),\eta)).$$

As $f_{1,\infty}$ is finitely generated, then set $$\{(f_{i}^{l})^{-1}(S(f_{i}^{l}(p_{x}),\eta))|i>0,n_{q}\geq l\geq 0\}$$
is finite.

Then $V(p_{x})$ is an open set in X which contain $p_{x}$ and satisfies for any $i>0,n_{q}\geq l\geq 0$, $$f_{i}^{l}(V(p_{x}))\subset S(f_{i}^{l}(p_{x}),\eta).$$

Because of the transitivity,
there is a positive integer $k$ satisfying $f_{1}^{k}(S(x,\epsilon))\bigcap V(p_{x})\not=\emptyset$ which also means that there is a $z\in S(x,\epsilon)$ making $f_{1}^{k}(z)\in V(p_{x})$ happen.

Let positive integer $j$ satisfy $(j-1)n_{q}<k\leq jn_{q}$, then$0\leq jn_{q}-k<n_{q}$.

As $f_{1}^{k}(z)\in V(p_{x})$,
$$f_{1}^{jn_{q}}(z)=f_{k+1}^{jn_{q}-k}(f_{1}^{k}(z))\in f_{k+1}^{jn_{q}-k}(V(p_{x}))\subset S(f_{k+1}^{jn_{q}-k}(p_{x}),\eta).$$

Hence $$d(f_{k+1}^{jn_{q}-k}(p_{x}),f_{k+1}^{jn_{q}-k}(f_{1}^{k}(z))=d(f_{k+11}^{jn_{q}}(p_{x}),f_{1}^{jn_{q}}(z))<\eta\eqno(a).$$

And for any positive integer $i,m$,$$d(x,f_{i}^{m}(orb(p_{x})))\geq d(x,orb(p_{x}))>\delta\geq 8\eta,$$
then $$d(f_{k+1}^{jn_{q}-k}(p_{x}),x)>8\eta\eqno(b).$$

For the triangle inequality and (a),(b),we have
\begin{eqnarray*}
\ d(f_{1}^{jn_{q}}(q),f_{1}^{jn_{q}}(z))&=&d(q,f_{1}^{jn_{q}}(z))\\
&=&d(q,f_{k+1}^{jn_{q}-k}(f_{1}^{k}(z)))\\
&\geq &d(f_{k+1}^{jn_{q}-k}(f_{1}^{k}(z)),x)-d(q,x)\\
&\geq &d(f_{k+1}^{jn_{q}-k}(f_{1}^{k}(p_{x})),x)-d(f_{k+1}^{jn_{q}-k}(f_{1}^{k}(p_{x})),f_{k+1}^{jn_{q}-k}(f_{1}^{k}(z)))-d(q,x)\\
&\geq &(8\eta -\eta)-\eta\\
&=&6\eta>0.
\end{eqnarray*}

And also because of the triangle inequality,
$$d(f_{1}^{jn_{q}}(q),f_{1}^{jn_{q}}(x))+d(f_{1}^{jn_{q}}(x),f_{1}^{jn_{q}}(x))\geq d(f_{1}^{jn_{q}}(q),f_{1}^{jn_{q}}(z))\geq 6\eta.$$

Then we have $d(f_{1}^{jn_{q}}(q),f_{1}^{jn_{q}}(x))\geq 3\eta>\eta$ or $d(f_{1}^{jn_{q}}(z),f_{1}^{jn_{q}}(x))\geq 3\eta>\eta$, and both $q$ and $z$ are in $S(x,\epsilon)$.

In conclusion, for the arbitrariness of $x$ and $\epsilon$, we know $\eta$ is a sensitivity constant for $(X,f_{1,\infty})$.$\Box$

\bigskip

Especially, when a NADS degenerates into an ADS, the Theorem  \ref{th3.7} just is Theorem\ref{th3.6}.

\bigskip

\begin{theorem}
A NADS which satisfying the conditions in Theorem \ref{th3.7} is Devaney chaos.
\end{theorem}

\bigskip

Now, we will share some results about the relationships between transitivity and sensitivities in linear NADS, which also could be seem as a special example of theorem \ref{th3.7}.

\begin{theorem}\label{th3.5}
Let $(X,f_{1,\infty})$ be a linear NADS.
If it is topological transitive, it has sensitivity dependence on initial conditions.
\end{theorem}

\bigskip

The proof of theorem\ref{th3.5} is totaly same as it in linear ADS which can be found in book\cite{Grosse} theorem2.30.

\bigskip
In fact, the notion of collective sensitivity comes from Wang's idea in paper\cite{Wang}. In that paper, the author show us that an ADS is collective sensitive is equivalent to its induced hyperspace system is sensitive dependence on initial conditions. We also confirm it is true in NDAS in paper\cite{yang}. Furthermore, in linear ADS, Chen\cite{CC} also show us that transitivity implies collective sensitivity. In the following three theorems, we will prove that in a class of linear NADS.

\bigskip

\begin{lemma}\label{l1}
Assume $(X,d,f_{1,\infty})$ is a commutative linear NADS. If it is topological transitive, there exists $\eta>0$ satisfying: for any $\epsilon>0$,
there exists $z_{1},z_{2}\in \{x\in X|d(x,0)<\epsilon\}$ and positive integer $k$ ensuring $d(f_{1}^{k}(z_{1}),f_{1}^{k}(z_{2}))>\eta$ established.
\end{lemma}

{\bf Proof:}
Let $s_{1},s_{2}\in X,s_{1}\not=s_{2}$ and $\delta=d(s_{1},s_{2})>0$.

Now we will prove that $\frac{\delta}{2}$ is just the $\eta$ we need.

For any $\epsilon>0$, note
$G_{1}=\{x\in X|d(x,s_{1})<\frac{\delta}{8}\},G_{2}=\{x\in X|d(x,s_{2})<\frac{\delta}{8}\},W=\{x\in X|d(x,0)<\epsilon\}.$

Since $f_{1,\infty}$ is transitive, there exists a positive integer $m$ ensuring $f_{1}^{m}(G_{1})\bigcap G_{2}\not= \emptyset$.

Hence there exists $y\in G_{1}$ which satisfy $f_{1}^{m}(y)\in G_{2}$.

As every $f_{i}$ is continues,
so is $f_{1}^{m}$. Then there is a $V(y)$ (which is a open neighborhood of $y$) satisfying
$f_{1}^{m}(V(y))\subset G_{2}$.

Let $G_{0}=G_{1}\bigcap V(y)$.
$G_{0}$ is an open set contained by
$G_{1}$
which satisfies\\
$$f_{1}^{m}(G_{0})\subset f_{1}^{m}(V(y))\subset G_{2}.$$
Because $(X,d,f_{1,\infty})$ is linear,
$f_{1}^{m}(0)=0\in W$.

Considering the continuity of $f_{1}^{m}$, there exists an open neighborhood of
0 : $U(0)\subset W$ ensuring
$f_{1}^{m}(U(0))\subset W$.

Also for the transitivity of $f_{1,\infty}$,
we can find a positive integer $k$ satisfying $f_{1}^{k}(U(0))\bigcap G_{0}\not= \emptyset$.

Therefore, there exists $y_{0}\in U(0)$ which can ensure that
$f_{1}^{k}(y_{0})\in G_{0}$.

So finally we have $f_{1}^{m}(f_{1}^{k}(y_{0}))\in f_{1}^{m}(G_{0})\subset G_{2}$.

For the commutativity of $f_{1,\infty}$, we have\\
$$f_{1}^{k}(f_{1}^{m}(y_{0}))=f_{1}^{m}(f_{1}^{k}(y_{0}))\in f_{1}^{m}(G_{0})\subset G_{2}.$$

Now we note $f_{1}^{m}(y_{0})$ and $y_{0}$ as $z_{1},z_{2}$.
Since $f_{1}^{m}(y_{0})\in f_{1}^{m}(U(0))\subset W,y_{0}\in U(0)\subset W$,
$z_{1},z_{2}\in W$ and
$$f_{1}^{k}(z_{1}))=f_{1}^{k}(f_{1}^{m}(y_{0}))\in f_{1}^{m}(G_{0})\subset G_{2},$$ $$f_{1}^{k}(z_{2})=f_{1}^{k}(y_{0})\in G_{0}\subset G_{1},$$

we have $d(z_{1},0)<\epsilon, d(z_{2},0)<\epsilon,d(f_{1}^{k}(z_{1}),f_{1}^{k}(z_{2})>\frac{\delta}{2}$.

Because of the arbitrariness of $\epsilon$, $\frac{\delta}{2}$ is just the $\eta$ we need.$\Box$

\bigskip

\begin{theorem}\label{3.5}
Assume $(X,d,f_{1,\infty})$ is a commutative linear NADS. If it is topological transitive, then it has collective sensitivity.
\end{theorem}

{\bf Proof:}
For any finitely many distinct points $x_{1},x_{2},x_{3},\cdots,x_{n}$ of $X$ and a arbitrary $\epsilon>0$, considering lemma\ref{l1}, we know there exists $\eta>0$ satisfying:
for $\epsilon$,
there exists $z_{1},z_{2}\in \{x\in X|d(x,0)<\epsilon\}$ and positive integer $k$ making

$$d(f_{1}^{k}(z_{1}),f_{1}^{k}(z_{2}))>\eta$$ established.

Now we are going to prove $\frac{\eta}{2}$ is a collective sensitivity constant of $(X,d,f_{1,\infty})$.

Since $z_{1},z_{2}\in \{x\in X|d(x,0)<\epsilon\}$,
$d(0,z_{1})<\epsilon,d(0,z_{2})<\epsilon.$

Let $y_{i}=x_{i}+z_{1},y_{i}^{*}=x_{i}+z_{2},i=1,2,3,\cdots,n$\\
then for any $i\in \{1,2,3,\cdots,n\}$,
$$d(x_{i},y_{i})=d(x_{i}-x_{i},y_{i}-x_{i})=d(0,z_{1})<\epsilon,$$
$$d(x_{i},y_{i}^{*})=d(x_{i}-x_{i},y_{i}^{*}-x_{i})=d(0,z_{2})<\epsilon.$$

As $d(f_{1}^{k}(z_{1}),f_{1}^{k}(z_{2}))>\eta,$
\begin{eqnarray*}
\ d(f_{1}^{k}(y_{i}),f_{1}^{k}(y_{i}^{*}))&=&d(f_{1}^{k}(x_{i}+z_{1}),f_{1}^{k}(x_{i}+z_{2}))\\
&=&d(f_{1}^{k}(x_{i})+f_{1}^{k}(z_{1}),f_{1}^{k}(x_{i})+f_{1}^{k}(z_{2}))\\
&=&d(f_{1}^{k}(z_{1}),f_{1}^{k}(z_{2}))>\eta,i=1,2,3,\cdots,n.
\end{eqnarray*}

For the triangle inequality, taking any given $x_{i_{0}}\in \{x_{1},x_{2},x_{3},\cdots ,x_{n}\}$, we always have $d(f_{1}^{k}(y_{i}),f_{1}^{k}(x_{i_{0}}))>\frac{\eta}{2}$ or $d(f_{1}^{k}(y_{i}^{*}),f_{1}^{k}(x_{i_{0}}))>\frac{\eta}{2}$ established.

That means there exists $x_{i}^{*}\in \{y_{i}^{*},y_{i}\},i=1,2,3,\cdots,n$ satisfying:
$$d(x_{i_{0}},x_{i}^{*})<\epsilon,$$
$$d(f_{1}^{k}(x_{i_{0}}),f_{1}^{k}(x_{i}^{*}))>\frac{\eta}{2}.$$

Hence $\frac{\eta}{2}$ is a collective sensitivity constant of $(X,d,f_{1,\infty})$.$\Box$

\bigskip

\begin{theorem}\label{syn}
Assume $(X,d,f_{1,\infty})$ is a commutative linear NADS. If it is topological transitive, then it has synchronous sensitivity.
\end{theorem}

{\bf Proof:}
For any finitely many distinct points $x_{1},x_{2},x_{3},\cdots,x_{n}$ of $X$ and a arbitrary $\epsilon>0$, for lemma \ref{l1}, we know there exists $\eta>0$ satisfying:
for $\epsilon$,
there exists $z_{1},z_{2}\in \{x\in X|d(x,0)<\epsilon\}$ and positive integer $k$ making
$$d(f_{1}^{k}(z_{1}),f_{1}^{k}(z_{2}))>\eta$$ established.

Now we will prove $\frac{\eta}{2}$ is a synchronous sensitivity constant of $(X,d,f_{1,\infty})$.

Just like theorem \ref{3.5}, let $y_{i}=x_{i}+z_{1},y_{i}^{*}=x_{i}+z_{2},i=1,2,3,\cdots,n$,
we have
$$d(x_{i},y_{i})=d(x_{i}-x_{i},y_{i}-x_{i})=d(0,z_{1})<\epsilon,$$
$$d(x_{i},y_{i}^{*})=d(x_{i}-x_{i},y_{i}^{*}-x_{i})=d(0,z_{2})<\epsilon,$$
\begin{eqnarray*}
\ d(f_{1}^{k}(y_{i}),f_{1}^{k}(y_{i}^{*}))&=&d(f_{1}^{k}(x_{i}+z_{1}),f_{1}^{k}(x_{i}+z_{2}))\\
&=&d(f_{1}^{k}(x_{i})+f_{1}^{k}(z_{1}),f_{1}^{k}(x_{i})+f_{1}^{k}(z_{2}))\\
&=&d(f_{1}^{k}(z_{1}),f_{1}^{k}(z_{2}))>\eta,i=1,2,3,\cdots,n.
\end{eqnarray*}
Because of the triangle inequality,
for any $x_{i},i\in\{1,2,3,\cdots,n\}$, we always have $d(f_{1}^{k}(y_{i}),f_{1}^{k}(x_{i}))>\frac{\eta}{2}$ or $d(f_{1}^{k}(y_{i}^{*}),f_{1}^{k}(x_{i}))>\frac{\eta}{2}$ established.

It also means there is $x_{i}^{*}\in \{y_{i}^{*},y_{i}\},i=1,2,3,\cdots,n$ satisfying:
$$d(x_{i},x_{i}^{*})<\epsilon,$$
$$d(f_{1}^{k}(x_{i}),f_{1}^{k}(x_{i}^{*}))>\frac{\eta}{2}.$$

That means $\frac{\eta}{2}$ is a synchronous sensitivity constant of $(X,d,f_{1,\infty})$.
$\Box$

\bigskip
The idea of synchronous sensitivity and the result of theorem \ref{syn} comes from the process of proving theorem\ref{3.5}.
Also, it is quite different with other sensitivities and stronger than sensitivity dependence on initial conditions.


\begin{thebibliography}{99}

\bibitem{Akin} E. Akin, S. Kolyada. Li-Yorke sensitivity. Nonlinearity, 2003; 16: 1421--1433.
\bibitem{Balibrea} F. Balibrea, P. Oprochab. Weak mixing and chaos in nonautonomous discrete systems. Applied Mathematics Letters, 2012; 25: 1135--1141.
\bibitem{Banks} J. Banks, J. Brooks, G. Gairns, G. Davis, D. Stacey. On Devaney's definition of chaos. The Amer. Math. Moonthly, 1992; 99: 332--334.
\bibitem{C2} J. S. Canovas, A. Linerob. Periodic structure of alternating continuous interval maps. Journal of difference Equations and Applications, 2006; 12(8): 847--858.
\bibitem{C1} J.S. Canovas. Li-Yorke chaos in a class of nonautonomous discrete systems. Journal of Difference Equations and Applications, 2011; 17(4):  479--486.
\bibitem{CC} C. Chen, Y. Wang. Some resarching about collective sensitivity in dynamical system. Northwest Universuty, 2013.
\bibitem{Devaney} R. L. Devaney. An introduction to chaotic dynamical systems. Addison-Wesley Publishing Company, 1989.
\bibitem{Dvorakova} J. Dvorakova. Chaos in nonautonomous discrete dynamical systems. Commun Nonlinear Sci Numer Simulat, 2012; 17: 4649--4652.
\bibitem{Glasner}  E. Glasner, B. Weiss, Sensitive dependence on initial conditions. Nonlinearity, 1993; 6: 1067--1075.
\bibitem{Grosse} K. Grosse-Erdmann, A. Peris. Linear Chaos. London: Springer-Verlag London, 2011.
\bibitem{Kolyada} S. Kolyada, L. Snoha. Topological entropy of non-autonomous dynamical system. Random Comp. Dyn., 1996; 4(2-3): 205--233.
\bibitem{Murillo} M. Murillo-Arcila, A. Peris. Mixing properties for nonautonomous linear dynamics and invariant sets. Applied Mathematics Letters, 2013; 26: 215--218.
\bibitem{Ruelle} D. Ruelle, F. Takens. On the nature of turbulence. Communications in Mathematical Physics, 1971; 20: 178--188.
\bibitem{Shi} Y. Shi, G. Chen. Chaos of time-varying discrete dynamical system. Journal of Difference Equations and Applicationns, 2009; 5: 429--449.
\bibitem{Silverman} S. Silverman, On maps with dense orbits and the definition of chaos. Rocky Mountain J. Math. 1992; 22: 353–-375.
\bibitem{Tian} C. Tian, G. Chen. Chaos of a sequence of maps in metric space. Chaos, Solitons and Fractals, 2006; 28: 1067--1075.
\bibitem{Wang}  Y. Wang, G. Wei, W. H. Campbell. Sensitive dependence on initial conditions between dynamical systems and their induced hyperspace dynamical system. Topology and its Applications, 2009; 156: 803--811.
\bibitem{Wu} X. Wu, P. Zhu. Chaos in a class of non-autonomous discrete system. Applied Mathematics Letters, 2013; 26: 431--426.
\bibitem{Xiong} J. Xiong. Chaos in the topological transitive systems. Science in China, 2005; 48: 929--939.
\bibitem{yang} C. Yang, Y. Wang, Z. Li, Some results about sensitivity on non-autonomous dynamical system in hyperspace. Pure and Applied Mathematics, 2014; 30(2): 201-206.
\end{thebibliography}
\end{document}